\documentclass[12pt]{amsart}
\usepackage[colorlinks=true,citecolor=black,linkcolor=black,urlcolor=blue]{hyperref}
\usepackage{amsmath}
\usepackage[english,  activeacute]{babel}
\usepackage[utf8]{inputenc}
\usepackage{amssymb}
\usepackage{amsthm}
\usepackage{graphics,graphicx}
\usepackage{enumerate}
\usepackage{array}
\usepackage{bm}
\usepackage{a4wide}
\usepackage{color, url}
\usepackage{float}
\usepackage{comment}
\usepackage{xcolor}
\allowdisplaybreaks

\usepackage{mmacells}
\usepackage{mma} 

\newcommand{\seqnum}[1]{\href{http://oeis.org/#1}{\underline{#1}}}
\DeclareMathOperator {\myRD}{RD}
\theoremstyle{plain}
\newtheorem{theorem}{Theorem}[section]

\newtheorem{lemma}[theorem]{Lemma}

\theoremstyle{definition}

\def\modd#1 #2{#1\ \mbox{\rm (mod}\ #2\mbox{\rm )}}

\usepackage[usestackEOL]{stackengine}[2013-10-15]
\stackMath

\title[]{On the arithmetic properties of partitions into parts simultaneously $4$-regular and $9$-distinct}

\date{\today}

\begin{document}

\author[M. L. Nadji]{Mohammed L. Nadji}
\address{Faculty of Mathematics, University of Science and Technology Houari Boumediene, RECITS laboratory, Algiers. LMAM laboratory, Jijel, 
Algeria}
\email{m.nadji@usthb.dz}

\author[M. Ahmia]{Moussa Ahmia}
\address{University of Mohamed Seddik Benyahia, LMAM laboratory, Jijel, 
Algeria}
\email{moussa.ahmia@univ-jijel.dz}

\begin{abstract}
In 2017, Keith presented a comprehensive survey on integer partitions into parts that are simultaneously regular, distinct, and/or flat. Recently, the authors initiated a study of partitions into parts that are simultaneously regular and distinct, examining them from both arithmetic and combinatorial perspectives. In particular, several Ramanujan-like congruences were obtained for $\myRD^{(\ell, t)}(n)$, the number of partitions of $n$ into parts that are simultaneously $\ell$-regular and $t$-distinct (parts appearing fewer than $t$ times), for various pairs $(\ell, t)$. In this paper, we focus on the case $(\ell, t)=(4,9)$ and conduct a thorough investigation of the arithmetic properties of $\myRD^{(4, 9)}(n)$. We establish several infinite families of congruences modulo $4$, $6$, and $12$, along with a collection of Ramanujan-like congruences modulo $24$.
\end{abstract}

\maketitle

\section{Introduction}
A \textit{partition} of a nonnegative integer $n$ is a finite sequence of positive integers which sum to $n$. We write $\lambda_{1}^{\alpha_{1}}+ \lambda_{2}^{\alpha_{2}}+\dots+ \lambda_{k}^{\alpha_{k}}=n$, where each part $\lambda_{i}$ occurs $\alpha_{i}$ times with $\lambda_{i}, \alpha_{i} \geq 1$ for $1 \leq i \leq k$. Here and throughout we will use the standard $q$-series notation
\begin{displaymath}
(a;q)_{n}:=\begin{cases} 
      \prod\limits_{i=0}^{n-1}(1-aq^{i}) &  \text{if $n>0$,} \\
      1 & \text{if $n=0$.}
    \end{cases}
\end{displaymath}
Moreover, $(a;q)_{\infty}=\lim_{n\rightarrow \infty}(a;q)_{n}$ for $|q|<1$, 
and $f_{k}$ is defined by 
\[
f_{k}:=(q^{k},q^{k})_{\infty}=\prod\limits_{n\geq 1}(1-q^{nk}).
\]

For a positive integer $\ell>1$, a partition is called $\ell$-regular if none of the parts are divisible by $\ell$. Let $b_{\ell}(n)$ denote the number of $\ell$-regular partitions of $n$. The generating function for $b_{\ell}(n)$ is given by
\[
\sum\limits_{n\geq 0} b_{\ell}(n) q^{n}= \dfrac{f_{\ell}}{f_{1}}.
\]
Regular partitions have been studied by many authors, including \cite{bb, r1, r2, r3, r4, r5, r6}.

An \textit{$\ell$-distinct partition} of  $n$ is a finite sequence of positive integers such that $\lambda_{1}^{\alpha_{1}}+ \lambda_{2}^{\alpha_{2}}+\dots+ \lambda_{k}^{\alpha_{k}}=n$, where $\lambda_{i} \in \mathbb{N}$, $1\leq \alpha_{i}< \ell$, and $\ell>2$. The number of $\ell$-distinct partitions of $n$, denoted by $p_{\ell}(n)$, satisfies the same generating function for $b_{\ell}(n)$. For more details on the $\ell$-distinct partitions, we refer the readers to \cite{r8, r9}.

A partition $\lambda=(\lambda_{1}^{\alpha_{1}},\lambda_{2}^{\alpha_{2}}, \dots, \lambda_{k}^{\alpha_{k}})$ of $n$ is said to be simultaneously $\ell$-regular and $t$-distinct if its parts $\lambda_{i}$ are indivisible by $\ell$ and $1\leq \alpha_{i} < t$ for $1\leq i \leq k$.  The number of such partitions of $n$ is denoted by $\myRD^{(\ell, t)}(n)$ and satisfies the generating function
\begin{equation} 
\sum_{n \geq 0} \myRD^{(\ell,t)}(n) q^{n} = \dfrac{f_{t}f_{\ell}}{f_{1}f_{\ell t}}. \label{mainequation} 
\end{equation}
For example, the nine partitions of $6$ into parts simultaneously $4$-regular and $9$-distinct are
\[
(6), (5,1), (3^2), (3,2,1), (3,1^3), (2^3), (2^2,1^2), (2,1^2), (1^6).
\]
This class of partitions was introduced by Keith \cite{r10}, who provided a comprehensive survey on partitions into parts that are simultaneously regular, distinct, and/or flat. The concept was later explored both combinatorially and arithmetically by the authors \cite{NadjiAhmia2024}. The number $\myRD^{(\ell,t)}(n)$ has been studied arithmetically for various pairs $(\ell, t)$ in \cite{NadjiAhmia2024}, including $(\ell, t)=(4,9)$, where several Ramanujan-like congruences  were established. For instance, the authors \cite{NadjiAhmia2024} presented the following congruences.

\begin{lemma}[Nadji and Ahmia \cite{NadjiAhmia2024}] For all $n \geq 0$, the following congruences hold \begin{align*} 
& \myRD^{(4, 9)}(4n+3) \equiv 0 \ (\mathrm{mod} \ 3),\\
& \myRD^{(4, 9)}(6n+2) \equiv 0 \ (\mathrm{mod} \ 2),\\
& \myRD^{(4, 9)}(6n+3) \equiv 0 \ (\mathrm{mod} \ 3),\\ 
& \myRD^{(4, 9)}(6n+4) \equiv 0 \ (\mathrm{mod} \ 4),\\ 
& \myRD^{(4, 9)}(6n+5) \equiv 0 \ (\mathrm{mod} \ 6),\\
& \myRD^{(4, 9)}(6n+7) \equiv 0 \ (\mathrm{mod} \ 12),\\
& \myRD^{(4, 9)}(8n+5) \equiv 0 \ (\mathrm{mod} \ 6). 
\end{align*} \end{lemma}
It is noteworthy that the enumeration function $\myRD^{(4, 9)}(n)$ corresponds to sequence \seqnum{A187020} (shifted by one term) in the Online Encyclopedia of Integer Sequences \cite{OEIS}.

In this paper, we delve into the arithmetic behavior of partitions into parts simultaneously $4$-regular and $9$-distinct by establishing some families of congruences for $\myRD^{(4, 9)}(n)$ modulo $4, 6$, and $12$, along with additional Ramanujan-like congruences modulo $24$.

\section{Preliminary results}
\label{Preres}
In this section, we list few dissection formulas which are useful in proving our main results. \emph{Ramanujan's general theta-function} $f(a,b)$ \cite[p. 34, 18.1]{Berndt1991} is defined by
\[
f(a,b):=\sum_{n=-\infty}^{\infty}a^{n(n+1)/2}b^{n(n-1)/2}=(-a,ab)_{\infty}(-b,ab)_{\infty}(ab,ab)_{\infty}.
\]
Some special cases of theta-function are denoted by
\[
\psi(q):=f(q,q^{3})=\sum_{n=0}^{\infty}q^{n(n+1)/2}=\dfrac{(q^{2},q^{2})_{\infty}}{(q,q^{2})_{\infty}}=\dfrac{f_{2}^{2}}{f_{1}},
\]
and
\[
f(-q):=f(-q^{2},-q)=\sum_{n=-\infty}^{\infty}(-1)^{n}q^{n(3n+1)/2}=(q;q)_{\infty}=f_{1}.
\]

\begin{lemma} \label{lemma1} The following $2$-dissections hold
\begin{align}
&\frac{f_{3}^{2}}{f_{1}^{2}}=\frac{f_{4}^{4}f_{6}f_{12}^{2}}{f_{2}^{5}f_{8}f_{24}}+2q\frac{f_{4}f_{6}^{2}f_{8}f_{24}}{f_{2}^{4}f_{12}}, \label{lemma1.1}\\
&\frac{f_{3}^{3}}{f_{1}}=\frac{f_{4}^{3}f_{6}^{2}}{f_{2}^{2}f_{12}}+q\frac{f_{12}^{3}}{f_{4}}. \label{lemma1.2}
\end{align}
\end{lemma}
Xia and Yao \cite{XY2} proved Equation \eqref{lemma1.1} and Hirschhorn et al. \cite{HGB} proved Equation \eqref{lemma1.2}.

The following result appears in the papers of Ramanujan   \cite[p. 212]{RAM}. 
\begin{lemma} \label{lemma2} We have the following $5$-dissection
    \begin{equation}
        f_{1}=f_{25}(a-q-q^{2}/a), \label{lemma2.1}
    \end{equation}
    where 
    $$a:= (q^{10}, q^{15}; q^{25})_{\infty}/(q^{10}, q^{5}; q^{20})_{\infty}.$$
\end{lemma}

\begin{lemma} \label{psidissectionlemma} For any odd prime $p$, we have
\begin{equation}
\psi(q)= \sum_{k=0}^{(p-3)/2} q^{k(k+1)/2} f(q^{\frac{p^{2}+(2k+1)p}{2}}, q^{\frac{p^{2}-(2k+1)p}{2}})+ q^{\frac{p^{2}-1}{8}}\psi(q^{p^{2}}).
 \label{psidissection}
\end{equation}
Furthermore, $\frac{m^{2}+m}{2}\not\equiv\frac{p^{2}-1}{8} \ (\mathrm{mod} \ p)$ for $0 \leq m \leq (p-3)/2$.
\end{lemma}

\begin{lemma} \label{fdissectionlemma} For any prime $p\geq 5$, we have
\begin{equation}
f_{1}= \sum_{\substack{k = (1-p)/2 \\ k\neq \frac{\pm p-1}{6}} }^{(p-1)/2} (-1)^{k} q^{k(3k+1)/2} f(-q^{\frac{3p^{2}+(6k+1)p}{2}}, -q^{\frac{3p^{2}-(6k+1)p}{2}})+(-1)^{\frac{\pm p-1 }{6}} q^{\frac{p^{2}-1}{24}}f_{p^{2}},
 \label{fdissection}
\end{equation}
where 
\[
   	\dfrac{\pm p-1}{6}=
    \begin{cases}
        \dfrac{p-1}{6}, \ if &  p \equiv \ 1 \ (\mathrm{mod} \ 6), \\
      	\dfrac{-p-1}{6}, \ if &  p \equiv \ -1 \ (\mathrm{mod}\ 6).
    \end{cases}
\]
\end{lemma}
Lemmas $\ref{psidissectionlemma}$ and $\ref{fdissectionlemma}$ are due to Cui and Gu \cite[Thm. 2.1 and Thm. 2.2]{r1}.

 \begin{lemma}\label{lemacon}
For all primes $p$ and all $k, m\geq 1$, we have 
\begin{equation}
 f_{pm}^{p^{k-1}} \equiv f_{m}^{p^{k}}  \ (\mathrm{mod} \ p^{k}). \label{binomiallemma} 
\end{equation}
\end{lemma}

Let $p$ be any odd prime and $\delta$ be any integer relatively prime to $p$. The \emph{Legendre symbol} $\bigl( \frac{\delta}{p} \bigl)$ is defined by
\begin{equation*}
        \left( \frac{\delta}{p} \right) = \begin{cases}
                      1,  &\text{if $\delta$ is a quadratic residue of $p$,}\\
                     -1,  &\text{if $\delta$ is a quadratic non-residue of $p$.}\\
                    \end{cases}
\end{equation*}

The following generating function dissections were proved by the authors \cite{NadjiAhmia2024}.
 \begin{lemma} \label{lemmanadjiahmia} We have
\begin{align}
&\sum_{n\geq 0} \myRD^{(4, 9)}(2n)q^{n}=\dfrac{f_{6}^{7}f_{9}^{7}}{f_{3}^{9}f_{18}^{5}}+2q\dfrac{f_{6}^{6}f_{9}^{4}}{f_{3}^{8}f_{18}^{2}}+4q^{2}\dfrac{f_{6}^{5}f_{9}f_{18}}{f_{3}^{7}}, \label{dd1}\\
&\sum_{n\geq 0} \myRD^{(4, 9)}(6n+2)q^{n}=2\dfrac{f_{2}^{6}f_{3}^{4}}{f_{1}^{8}f_{6}^{2}}, \label{d1}\\
&\sum_{n\geq 0} \myRD^{(4, 9)}(6n+4)q^{n}=4\dfrac{f_{2}^{5}f_{3}f_{6}}{f_{1}^{7}}, \label{d2}\\
&\sum_{n\geq 0} \myRD^{(4, 9)}(6n+5)q^{n}= 6\dfrac{f_{2}^{2}f_{3}^{2}f_{6}^{2}}{f_{1}^{6}}. \label{d3}
\end{align}
\end{lemma}

\section{Congruences for $\myRD^{(4, 9)}(n)$ modulo $4$}

\begin{theorem}\label{teor43g} For any prime $p \equiv 3 \ (\mathrm{mod} \ 4)$, $\alpha\geq 0$, and $n\geq 0$, we have
\begin{equation*}
\myRD^{(4, 9)}\Bigl( 12p^{2\alpha+1}(pn+i)+ p^{2\alpha+2}+1 \Bigl) \equiv 0 \ (\mathrm{mod} \ 4),
\end{equation*}
for all  $1\leq i \leq p-1$.
\end{theorem}

\begin{proof} 
In view of \eqref{binomiallemma} and \eqref{dd1}, with $p=2$ and $k=1$, we have
\[
    \sum_{n\geq 0} \myRD^{(4, 9)}(2n) q^{n}=\dfrac{f_{6}^{7}f_{9}^{7}}{f_{3}^{9}f_{18}^{5}}+2q\dfrac{f_{6}^{6}f_{9}^{4}}{f_{3}^{8}f_{18}^{2}}+4q^{2}\dfrac{f_{6}^{5}f_{9}f_{18}}{f_{3}^{7}} \equiv \dfrac{f_{6}^{7}f_{9}^{7}}{f_{3}^{9}f_{18}^{5}}+2qf_{6}^{2} \ (\mathrm{mod} \ 4).
\]
If we extract the terms of the form $q^{3n+1}$ from both sides of the above equation, we get
\[
\sum_{n\geq 0} \myRD^{(4, 9)}(6n+2) q^{n} \equiv 2f_{2}^{2} \ (\mathrm{mod} \ 4).
\]
That is,
\begin{equation}
    \sum_{n\geq 0} \myRD^{(4, 9)}(12n+2) q^{n} \equiv 2f_{1}^{2} \ (\mathrm{mod} \ 4). \label{XX1}
\end{equation}
Define
\begin{equation} 
\sum_{n\geq 0} a(n)q^{n} = f_{1}^{2}. \label{c1}
\end{equation}
Combining Equations \eqref{XX1} and \eqref{c1}, we deduce that
\begin{equation}
\myRD^{(4, 9)}(12n+2) \equiv  2a(n) \ (\mathrm{mod} \ 4). \label{c2}
\end{equation}
Now, we consider the following congruence equation:
\[
\frac{3k^{2}+k}{2}+\frac{3m^{2}+m}{2} \equiv \frac{p^{2}-1}{12} \ (\mathrm{mod} \ p).
\]
It is equivalent to
\begin{equation}
(6k+1)^{2}+(6m+1)^{2} \equiv 0 \ (\mathrm{mod} \ p), \label{c3}
\end{equation}
where $-(p-1)/2 \leq k,m \leq (p-1)/2$ and $p$ is a prime such that $(\frac{-1}{p})=-1$. Since $(\frac{-1}{p})=-1$ for $p\equiv 3 \ (\mathrm{mod} \ 4)$, the congruence relation \eqref{c3} holds if and only if both $k=m=(\pm p-1)/6$. Substituting Equation \eqref{fdissection} into \eqref{c1} and then extracting the terms in which the powers of $q$ are congruent to $(p^{2}-1)/12$ modulo $p$, and then divide by $q^{\frac{p^{2}-1}{12}}$, we find that
\[
\sum_{n\geq 0} a\left(pn+\frac{p^{2}-1}{12}\right)q^{pn}=f_{p^{2}}^{2}.
\]
Therefore,
\begin{equation*}
\sum_{n\geq 0} a\left(p^{2}n+\frac{p^{2}-1}{12}\right) q^{n} =f_{1}^{2},  
\end{equation*}
and for all $n\geq 0$,
\begin{equation}
a \left(p^{2}n+pi+\frac{p^{2}-1}{12}\right)=0, \label{c5}
\end{equation}
where  $1\leq i \leq p-1$. By induction we can prove that for $n\geq 0$ and $\alpha\geq 0$,
\begin{equation}
a\left(p^{2\alpha}n+\frac{p^{2\alpha}-1}{12}\right)=a(n). \label{c6}
\end{equation}
Replacing $n$ by $p^{2}n+pi+\frac{p^{2}-1}{12}$ ($1\leq i \leq p-1$) in Equation \eqref{c6} and using \eqref{c5}, we find that for all nonnegative integers $n$ and $\alpha$,
\[
a \left(p^{2\alpha+2}n+p^{2\alpha+1}i+\frac{p^{2\alpha+2}-1}{12}\right)=0. 
\]
Finally, replacing $n$ by $p^{2\alpha+2}n+p^{2\alpha+1}i+\frac{p^{2\alpha+2}-1}{12}$ in Equation \eqref{c2} ($1\leq i \leq p-1$) we obtain the desired result.
\end{proof}

\section{Congruences for $\myRD^{(4, 9)}(n)$ modulo $6$}
 \begin{theorem} For all integers $\alpha \geq 0$ and $n\geq 0$,
 \begin{equation}
\myRD^{(4, 9)}(6\cdot 5^{2\alpha+2}n+6\cdot 5^{2\alpha+1}i+ 5^{2\alpha+2} +1 ) \equiv 0 \ (\mathrm{mod} \ 6), \label{famcongmod6}
\end{equation}
where $i=1,2,3,4.$
\end{theorem}

\begin{proof}
In view of \eqref{binomiallemma} and \eqref{d1}, with $p=3$ and $k=1$, we have
\begin{equation}
    \sum_{n\geq 0} \myRD^{(4, 9)}(6n+2) q^{n}=2\dfrac{f_{2}^{6}f_{3}^{4}}{f_{1}^{8}f_{6}^{2}} \equiv 2f_{1}^{4} \ (\mathrm{mod} \ 6). \label{aa1}
\end{equation}
Employing \eqref{lemma2.1} into \eqref{aa1}, and then extracting the terms of the form $q^{5n+4}$ from both sides of the resulting equation, we get
\begin{equation}
    \sum_{n\geq 0} \myRD^{(4, 9)}(30n+26) q^{n} \equiv 2f_{5}^{4} \ (\mathrm{mod} \ 6), \label{aa2}
\end{equation}
which implies that
\begin{equation}
    \sum_{n\geq 0} \myRD^{(4, 9)}(150n+26) q^{n} \equiv 2f_{1}^{4} \ (\mathrm{mod} \ 6). \label{aa3}
\end{equation}
From \eqref{aa1} and \eqref{aa3}, we find that
\[
    \myRD^{(4, 9)}(6n+2) \equiv \myRD^{(4, 9)}(150n+26)  \ (\mathrm{mod} \ 6). 
\]
Utilizing \eqref{aa3} and mathematical induction on $\alpha \geq 0$, we obtain
\begin{equation}
\myRD^{(4, 9)}(6n+2) \equiv \myRD^{(4, 9)}(6\cdot 5^{2\alpha+2}n+ 5^{2\alpha+2} +1)  \ (\mathrm{mod} \ 6).  \label{aa4}
\end{equation}
From \eqref{aa2}, we find that
\begin{equation}
     \myRD^{(4, 9)}(150n+30i+26)  \equiv 0 \ (\mathrm{mod} \ 6),  \ i=1,2,3,4. \label{aa5}
\end{equation}
Using Equations \eqref{aa4} and \eqref{aa5}, we obtain the desired result \eqref{famcongmod6}.
\end{proof}

 \section{Congruences for $\myRD^{(4, 9)}(n)$ modulo $12$}
 \begin{theorem} \label{thm} For any prime $p\equiv 5 \ (\mathrm{mod} \ 6)$, $\alpha \geq 0$, and $n\geq 0$, we have
\[
\myRD^{(4, 9)} \big(6p^{2\alpha+1}(pn+i)+3p^{2\alpha+2}+1 \big)\equiv 0 \ (\mathrm{mod} \ 12),
\]
for all $1\leq i \leq p-1$.
\end{theorem}

\begin{proof} 
In view of \eqref{binomiallemma} and \eqref{d2}, with $p=3$ and $k=1$, we find that
\begin{equation}
\sum_{n\geq 0} \myRD^{(4, 9)}(6n+4)q^{n}=4\dfrac{f_{2}^{5}f_{3}f_{6}}{f_{1}^{7}} \equiv 4 \dfrac{f_{2}^{2} f_{6}^{2}}{f_{1} f_{3}} = 4\psi(q) \psi(q^{3}) \ (\mathrm{mod} \ 12). \label{lel1}
\end{equation}
Define
\begin{equation} 
\sum_{n\geq 0} b(n)q^{n} = \psi(q) \psi(q^{3}). \label{e1}
\end{equation}
From  Equations \eqref{lel1} and \eqref{e1} we have 
\begin{equation}
\myRD^{(4, 9)}(6n+4) \equiv 4b(n) \ (\mathrm{mod} \ 12).  \label{e2}
\end{equation} 
Now, consider the congruence equation
\[
\frac{k^{2}+k}{2}+3\cdot \frac{m^{2}+m}{2} \equiv  \frac{4p^{2}-4}{8} \ (\mathrm{mod} \ p).
\]
which is equivalent to
\begin{equation}
(2k+1)^{2}+ 3\cdot (2m+1)^{2}  \equiv  0 \ (\mathrm{mod} \ p), \label{e3}
\end{equation}
where $0 \leq k,m \leq (p-1)/2$ and $p$ is a prime number such that ($\frac{-3}{p}$)$=-1$. Since ($\frac{-3}{p}$)$=-1$ for $p \equiv 5 \ (\mathrm{mod} \ 6)$, the congruence relation of Equation \eqref{e3} holds if and only if both $k=m=(p-1)/2$. Substitute Equation \eqref{psidissection} into \eqref{e3} and extract the terms in which the powers of $q$ are congruent to $(p^{2}-1)/2$ modulo $p$, and then divide by $q^{\frac{p^{2}-1}{2}}$, we find that
\[
\sum_{n\geq 0} b \left(pn+\frac{p^{2}-1}{2} \right) q^{pn} = \psi(q^{p^{2}}) \psi(q^{3p^{2}}),
\]
which implies that 
\[
\sum_{n\geq 0} b\left( p^{2}n+\frac{p^{2}-1}{2} \right)q^{n} = \psi(q) \psi(q^{3}),
\]
and for $n\geq 0$,
\begin{equation}
b \left(p^{2}n+pi+\frac{p^{2}-1}{2}\right) = 0, \label{e5}
\end{equation}
where  $1\leq i \leq p-1$. By induction, we obtain  that for all $n, \alpha \geq 0$,
\begin{equation}
b\left( p^{2\alpha}n+\frac{p^{2\alpha}-1}{2} \right)=b(n). \label{e6}
\end{equation}
Replacing $n$ by $p^{2}n+pi+\frac{p^{2}-1}{2}$ ($1\leq i \leq p-1$) in \eqref{e6} and using \eqref{e5}, we deduce that for $n\geq 0$ and $\alpha \geq 0$,
\[
b\left( p^{2\alpha+2}n+p^{2\alpha+1}i+\frac{p^{2\alpha+2}-1}{2} \right)=0.
\]
Replacing $n$ by $p^{2\alpha+2}n+p^{2\alpha+1}i+\frac{p^{2\alpha+2}-1}{2}$ in Equation \eqref{e2}, we obtain the desired result.
\end{proof}

\section{Congruences for $\myRD^{(4, 9)}(n)$ modulo $24$}
\begin{theorem} For all $n\geq 0$,
\begin{align}
& \myRD^{(4, 9)}(24n+23) \equiv  0 \ (\mathrm{mod} \ 24), \label{cong24-1}\\
& \myRD^{(4, 9)}(48n+29) \equiv  0 \ (\mathrm{mod} \ 24),\label{cong24-2}\\
& \myRD^{(4, 9)}(96n+89) \equiv 0 \ (\mathrm{mod} \ 24). \label{cong24-3}
\end{align}
\end{theorem}
 
 \begin{proof} In view of \eqref{binomiallemma} and \eqref{d3}, with $p=k=2$, we obtain
\[
    \sum_{n\geq 0} \myRD^{(4, 9)}(6n+5)q^{n}= 6\dfrac{f_{2}^{2}f_{3}^{2}f_{6}^{2}}{f_{1}^{6}} \equiv 6 \dfrac{f_{3}^{2} f_{6}^{2}}{f_{1}^{2}} \ (\mathrm{mod} \ 24).
\]
 Substituting \eqref{lemma1.1} into the above equation, we get
\begin{equation}
    \sum_{n\geq 0} \myRD^{(4, 9)}(6n+5)q^{n} \equiv 6 \dfrac{f_{4}^{4} f_{6}^{3} f_{12}^{2}}{f_{2}^{5} f_{8} f_{24}}+12q \dfrac{f_{4} f_{6}^{4} f_{8} f_{24}}{f_{2}^{4} f_{12}} \ (\mathrm{mod} \ 24). \label{m1}
\end{equation}
If we extract the even and the odd powers from both sides of \eqref{m1}, we find that
\begin{equation}
    \sum_{n\geq 0} \myRD^{(4, 9)}(12n+5)q^{n} \equiv 6 \dfrac{f_{2}^{4} f_{3}^{3} f_{6}^{2}}{f_{1}^{5} f_{4} f_{12}} \equiv 6 \dfrac{f_{2}^{2} f_{3}^{3} f_{6}^{2}}{f_{1} f_{4} f_{12}} \ (\mathrm{mod} \ 24), \label{m2}
\end{equation}
and
\begin{equation}
    \sum_{n\geq 0} \myRD^{(4, 9)}(12n+11)q^{n} \equiv 12 \dfrac{f_{2} f_{3}^{4} f_{4} f_{12}}{f_{1}^{4} f_{6}} \equiv 12 \dfrac{f_{4} f_{6} f_{12}}{f_{2}} \ (\mathrm{mod} \ 24). \label{m3}
\end{equation}
Employing \eqref{lemma1.2} in \eqref{m2}, we arrive at
\begin{equation}
    \sum_{n\geq 0} \myRD^{(4, 9)}(12n+5)q^{n} \equiv 6 \dfrac{f_{4}^{2} f_{6}^{4}}{f_{12}^{2}} + 6 q \dfrac{f_{2}^{2}  f_{6}^{2} f_{12}^{2}}{f_{4}^{2}} \ (\mathrm{mod} \ 24). \label{m4}
\end{equation}
Extracting the odd and the even powers from both sides of \eqref{m4}, we obtain
\begin{equation}
    \sum_{n\geq 0} \myRD^{(4, 9)}(24n+5)q^{n} \equiv 6 \dfrac{f_{2}^{2} f_{3}^{4}}{f_{6}^{2}} \equiv 6 f_{2}^{2} \ (\mathrm{mod} \ 24), \label{m5}
\end{equation}
and
\begin{equation}
    \sum_{n\geq 0} \myRD^{(4, 9)}(24n+17)q^{n} \equiv 6 \dfrac{f_{1}^{2}  f_{3}^{2} f_{6}^{2}}{f_{2}^{2}} \equiv 6 \dfrac{f_{3}^{2} f_{6}^{2}}{f_{1}^{2}} \ (\mathrm{mod} \ 24). \label{m6}
\end{equation}
Congruence \eqref{cong24-1} follows from \eqref{m3} by collecting the terms involving $q^{2n+1}$ from both sides, and congruence \eqref{cong24-2} follows from \eqref{m5} by equating the odd powers on both sides.

If we substitute \eqref{lemma1.1} into \eqref{m6}, and then extracting the odd powers from both sides of the resulting equation, we arrive at
\begin{equation}
    \sum_{n\geq 0} \myRD^{(4, 9)}(48n+41)q^{n} \equiv  12 \dfrac{f_{4} f_{6} f_{12}}{f_{2}} \ (\mathrm{mod} \ 24). \label{m7}
\end{equation}
Congruence \eqref{cong24-3} is an immediate result from \eqref{m7}.
 \end{proof}

\section*{Acknowledgment}
The authors would like to thank the reviewers for their valuable remarks and suggestions to improve the original manuscript. This work was supported by DG-RSDT (Algeria), PRFU Project, No. C00L03UN180120220002. 

\end{document}